\documentclass[11pt]{article}
\usepackage{graphicx}
\usepackage{amsfonts}
  \def\fH{{\cal H}}

\usepackage{theorem}
\newtheorem{Lem}{Lemma}[section]

\newtheorem{The}[Lem]{Theorem}
\newtheorem{Prop}[Lem]{Proposition}
\newtheorem{Cor}[Lem]{Corollary}

\newtheorem{Rem}[Lem]{Remark}

\newcommand{\qed}{\hbox{\rule{6pt}{6pt}}}

\begin{document}
\title{On refined Young inequalities and reverse inequalities}
\author{Shigeru Furuichi\footnote{E-mail:furuichi@chs.nihon-u.ac.jp} \\
{\small Department of Computer Science and System Analysis,}\\
{\small College of Humanities and Sciences, Nihon University,}\\
{\small 3-25-40, Sakurajyousui, Setagaya-ku, Tokyo, 156-8550, Japan}}
\date{}
\maketitle

{\bf Abstract.} In this paper, we show refined Young inequalities for two positive operators. 
Our results refine the ordering relations among the arithmetic mean, the geometric mean and the harmonic mean for two positive operators.
In addition, we give two different reverse inequalities for the refined Young inequality for two positive operators.
\vspace{3mm}

{\bf Keywords : } Young inequality, reverse inequality, positive operator, operator mean, Specht's ratio and operator inequality 
\vspace{3mm}

{\bf 2000 Mathematics Subject Classification : } 15A45
\vspace{3mm}

\maketitle


\section{Introduction}
We start from the famous arithmetic-geometric mean inequality which is often called Young inequality:
\begin{equation} \label{ineq01_rev}
(1-\nu) a+ \nu b \geq a^{1-\nu} b^{\nu}
\end{equation}
for nonnegative real numbers $a$, $b$ and $\nu \in [0,1]$.

Recently, the inequality (\ref{ineq01_rev}) was refined by F.Kittaneh and Y.Manasrah in the following form,
for the purpose of the study on matrix norm inequalities.
\begin{Prop} {\bf (\cite{BK,KM,IN})}    \label{prop02}
For $a, b\geq 0$ and $\nu \in [0,1]$, we have
\begin{equation} \label{ineq02_rev}
(1-\nu) a+ \nu b \geq a^{1-\nu} b^{\nu}+ r (\sqrt{a}-\sqrt{b})^2,
\end{equation}
where $r\equiv \min\{\nu,1-\nu\}$.
\end{Prop}
It is notable that the inequality (\ref{ineq02_rev}) was first proved in (6.32) on page 46 of the reference \cite{BK}.

In the section 2 of this paper,  we give refined Young inequalities for two positive operators based on the scalar inequality (\ref{ineq02_rev}).

As for the reverse inequalities of the Young inequality, M.Tominaga gave the following interesting operator inequalities. 
He called them {\it converse} inequalities, however we use the term {\it reverse} for such inequalities, throughout this paper.
\begin{Prop} {\bf (\cite{Tom1})}
Let $\nu \in [0,1]$, positive operators $A$ and $B$ such that $0<mI \leq  A,B \leq MI $ with $h\equiv \frac{M}{m}>1$.
Then we have the following inequalities for every  $\nu \in [0,1]$.
\begin{itemize}
\item[(i)] (Reverse ratio inequality)
$$
S(h)  A \sharp_{\nu}B \geq (1-\nu) A + \nu B,
$$
where the constant $S(h)$ is called Specht's ratio \cite{Specht,JIFUJII} and defined by
$$
S(h) \equiv \frac{h^{\frac{1}{h-1}}}{e\log{h^{\frac{1}{h-1}}}},\quad (h \neq 1)
$$
for positive real number $h$.
\item[(ii)] (Reverse difference inequality)
\begin{equation}\label{prop_rev_diff_ineq}
h  L(m,M)\log S(h) +A \sharp_{\nu}B \geq  (1-\nu) A + \nu B, 
\end{equation}
where the logarithmic mean $L$ is defined by
$$
L(x,y)\equiv \frac{y-x}{\log y-\log x},\,\,(x\neq y)\quad L(x,x)\equiv x
$$
for two positive real numbers $x$ and $y$.
\end{itemize}
\end{Prop}

In the section 3 of this paper, we give reverse ratio type inequalities of the refined Young inequality for positive operators.
In the section 4 of this paper, we also give reverse difference type inequalities of the refined Young inequality for positive operators.

\section{Refined Young inequalities for positive operators}
Let $\fH$ be a complex Hilbert space. We also represent the set of all bounded operators on $\fH$ by $B(\fH)$.
If $A\in B(\fH)$ satisfies $A^*=A$, then $A$ is called a self-adjoint operator.
A self-adjoint operator $A$ satisfies $\langle x \vert A \vert x\rangle \geq 0$ for any $\vert x \rangle \in \fH$, then $A$ is called a positive operator.
For two self-adjoint operators $A$ and $B$, $A\geq B$ means $A-B\geq 0$.

It is well-known that we have the following Young inequalities for invertible positive operators $A$ and $B$:
\begin{equation} \label{orig_Young_ineq}
(1-\nu) A+ \nu B \geq A\sharp_{\nu}B \geq \left\{(1-\nu) A^{-1}+ \nu B^{-1} \right\}^{-1},
\end{equation}
where $A\sharp_{\nu}B \equiv A^{1/2}(A^{-1/2}BA^{-1/2})^{\nu}A^{1/2}$  defined for $\nu\in[0,1]$.
The power mean was originally introduced in the paper \cite{KA}.
The simplified and elegant proof for the inequalities (\ref{orig_Young_ineq}) 
was given in \cite{FY}. See also \cite{Furuta} for the reader having interests in operator
inequalities.

As a refinement of the inequalities (\ref{orig_Young_ineq}),  we have the following refined Young inequality for positive operators.
\begin{The} \label{the01}
For $\nu \in [0,1]$ and positive operators $A$ and $B$, we have
\begin{eqnarray}  
 (1-\nu )A+\nu B  &\geq& A\sharp_{\nu} B  +  2r \left(\frac{A+B}{2} -A\sharp_{1/2}B \right) \label{the01-ineq01} \\
&\geq & A\sharp_{\nu} B  \label{the01-ineq02} \\
&\geq& \left\{ A^{-1}\sharp_{\nu} B^{-1}+2r \left(\frac{A^{-1}+B^{-1}}{2} -A^{-1}\sharp_{1/2}B^{-1} \right)\right\}^{-1} \label{the01-ineq03} \\
&\geq& \left\{  (1-\nu )A^{-1}+\nu B^{-1}   \right\}^{-1}   \label{the01-ineq04} 
\end{eqnarray}
where $r\equiv \min\left\{\nu,1-\nu\right\}$ and  $A\sharp_{\nu}B \equiv A^{1/2}(A^{-1/2}BA^{-1/2})^{\nu}A^{1/2}$  defined for $\nu\in[0,1]$.
\end{The}

To prove Theorem \ref{the01}, we use the following lemma. 
\begin{Lem} \label{lem_hg}
For invertible positive operators $X$ and $Y$, we have
$$
(X+Y)^{-1} = X^{-1} -X^{-1}(X^{-1}+Y^{-1})^{-1}X^{-1}.
$$
\end{Lem}
{\it Proof}:
Since $(X+Y)(X+Y)^{-1}X=X$, we have $  X(X+Y)^{-1}X  +  Y(X+Y)^{-1}X =X$.
Thus we have
$X(X+Y)^{-1}X =X-Y(X+Y)^{-1}X =X-(X^{-1}(X+Y)Y^{-1})^{-1}=X-(X^{-1}+Y^{-1})^{-1}$.
Multiplying $X^{-1}$ from both sides, we obtain the lemma.
\hfill \qed

{\it Proof of Theorem \ref{the01}}:
The second inequality (\ref{the01-ineq02}) is clear, since we have $2r \left(\frac{A+B}{2} -A\sharp_{1/2}B \right) \geq 0$. We prove the first inequality.
From the inequality (\ref{ineq02_rev}), we have for $\nu \in [0,1]$ and $x \geq 0$
$$\nu x+1-\nu-x^{\nu} -r(\sqrt{x}-1)^2 \geq 0.$$
By the standard operational calculus, we have
\begin{eqnarray}
\nu T +1-\nu &\geq& T^{\nu} +r (T^{1/2}-1)^2 \nonumber \\
&=&T^{\nu} +r (T-2T^{1/2}+1) \label{ineq01}
\end{eqnarray}
for a positive operator $T$ and  $\nu \in [0,1]$.
From here, we suppose that $A$ is an invertible.
(For a general case, we consider the invertible positive operator $A_{\epsilon}\equiv A+\epsilon I$
for positive real number $\epsilon$. If we take a limit $\epsilon \to 0$, the following result also holds.
Throughout this paper, we apply this continuity argument, however, from now on, we omit such descriptions for simplicity.)
Substituting $T=A^{-1/2}BA^{-1/2}  $ into the inequality (\ref{ineq01}), we have
$$
\nu A^{-1/2}BA^{-1/2}   +1-\nu \geq \left(A^{-1/2}BA^{-1/2}\right)  ^{\nu} +r \left\{ A^{-1/2}BA^{-1/2}   
-2\left(A^{-1/2}BA^{-1/2}\right)  ^{1/2}   +1 \right\}
$$
Multiplying $A^{1/2}$ to the above inequality from both sides, we have
$$
(1-\nu) A+ \nu B \geq A\sharp_{\nu}B +r \left( A+B -2 A \sharp_{1/2}B\right),
$$
which proves the inequality (\ref{the01-ineq01}).

Replacing  $A$ and $B$ by $A^{-1}$ and $B^{-1}$ in the inequality (\ref{the01-ineq01}), 
respectively and taking the inverse of both sides, then we have the last inequality (\ref{the01-ineq04}).

By Lemma \ref{lem_hg}, the right hand side of the inequality (\ref{the01-ineq03}) can be calculated as
$$
R.H.S.= A\sharp_{\nu}B -\left(A\sharp_{\nu}B\right) \left[ A\sharp_{\nu}B + \left\{ 2r\left( \frac{A^{-1}+B^{-1}}{2}-A^{-1}\sharp_{1/2}B^{-1}\right) \right\}^{-1} \right]^{-1} \left(A\sharp_{\nu}B\right). 
$$
Since, $\left(A^{-1}\sharp_{\nu}B^{-1}\right)^{-1}=A\sharp_{\nu}B  \geq 0$  
and $2r\left( \frac{A^{-1}+B^{-1}}{2}-A^{-1}\sharp_{1/2}B^{-1}\right)\geq 0$, we have the third inequality (\ref{the01-ineq03}),
which completes the proof.

\hfill \qed

In the paper \cite{Furuta1}, the equivalent relation between the Young inequality and the H\"older-McCarthy inequality \cite{Mc} 
was shown by a simplified elegent proof.
Here we show a kind of the refinement of the H\"older-McCarthy inequality applying Theorem \ref{the01}.

\begin{Cor}
For $\nu\in[0,1]$ and any  positive operator $A$ on the Hilbert space $\fH$ and any unit vector $\vert x \rangle \in \fH$, if
$\langle x \vert A \vert x \rangle \neq 0$, then we have
\begin{equation} \label{ineq02}
1- \langle x \vert A\vert x\rangle ^{-\nu}  \langle x \vert A^{\nu }\vert x\rangle
\geq r\left( 1- \langle x \vert A\vert x\rangle ^{-1/2}  \langle x \vert A^{1/2 }\vert x\rangle   \right)^2,
\end{equation}
where $r\equiv \min\left\{\nu,1-\nu\right\}$.
\end{Cor}
{\it Proof}:
If $\nu =0$, then the inequality (\ref{ineq02}) is trivial. It is sufficient that we prove it for the case of $\nu\in (0,1]$.
In the inequality (\ref{ineq01}) ,
we put $T=k^{\frac{1}{\nu}}A$, for any positive real number $k$ and by the unit vector $\vert x \rangle \in \fH$, we have
\begin{equation}\label{ineq03}
\nu k^{\frac{1}{\nu}}  \langle x \vert A \vert x \rangle +1-\nu \geq k  \langle x \vert A^{\nu} \vert x \rangle
+r\left( k^{\frac{1}{2\nu}}     \langle x \vert A^{1/2} \vert x \rangle  -1\right)^2
\end{equation}
In the inequality (\ref{ineq03}), if we put $k=\langle x \vert A \vert x \rangle ^{-\nu}$, then we obtain the inequality (\ref{ineq02}).

\hfill \qed

\begin{Rem}
From H\"older-McCarthy inequality \cite{Mc}:
\begin{equation}  \label{holder01}
\langle x \vert A \vert x \rangle ^{\nu} \geq \langle x\vert A^{\nu} \vert x\rangle
\end{equation}
for any unit vector $\vert x \rangle \in \fH$,
if $\langle x \vert A \vert x \rangle \neq 0$, then we have
$$
1- \langle x \vert A\vert x\rangle ^{-\nu}  \langle x \vert A^{\nu }\vert x\rangle  \geq 0.
$$
The inequality (\ref{ineq02}) gives a refined one for the above inequality 
which is equivalent to the inequality (\ref{holder01}) in the case of $\langle x \vert A \vert x \rangle \neq 0$.
\end{Rem}


\section{A reverse ratio inequality for a refined Young inequality}
For positive real numbers $a,b$ and $\nu\in[0,1]$, M.Tominaga showed the following inequality \cite{Tom1}:
\begin{equation}\label{ineq21} 
S\left(\frac{a}{b} \right) a^{1-\nu} b^{\nu} \geq (1-\nu) a + \nu b,
\end{equation}
which is called the converse ratio inequality for the Young inequality in \cite{Tom1}.
In this section, we show the reverse ratio inequality of the refined Young inequality (\ref{ineq02_rev}).
\begin{Lem}   \label{lem21}
For positive real numbers $a, b$ and $\nu \in [0,1]$, we have
\begin{equation}\label{ineq22}
S\left(\sqrt{\frac{a}{b}} \right)a^{1-\nu}b^{\nu} \geq (1-\nu) a + \nu b -r(\sqrt{a}-\sqrt{b})^2,
\end{equation}
where $r\equiv \min\left\{\nu,1-\nu \right\}$.
\end{Lem}
{\it Proof}:
\begin{itemize}
\item[(i)]For the case of $\nu \leq1/2$, $r=\nu$.
We consider the following function.
$$
g_b(\nu) \equiv \frac{\nu b+(1-\nu) -\nu(\sqrt{b} -1)^2}{b^{\nu}},\quad \left(0\leq \nu \leq \frac{1}{2}\right) .
$$
Then we have 
$$
g'_b(\nu) = \frac{2(\sqrt{b}-1)-\left\{2(\sqrt{b}-1)\nu + 1\right\}\log b}{b^{\nu}}
$$
so that the equation $g'_b(\nu) = 0$ implies 
$$\nu =  \frac{1}{\log b}-\frac{1}{2(\sqrt{b}-1)} \equiv \nu_b.$$
From the Klein inequality:
$$
1-\frac{1}{\sqrt{b}} \leq \log \sqrt{b} \leq \sqrt{b} -1,\quad (b>0)
$$
we have $\nu_b \in [0,\frac{1}{2}]$.
We also find that $g'_b(\nu) >0$ for $\nu < \nu_b$ and $g'_b(\nu) <0$ for $\nu > \nu_b$.
Thus the function $g_b(\nu)$ takes a maximum value when  $\nu=\nu_b,\, (b\neq 1)$ and it is calculated as follows.
\begin{eqnarray*}
\max_{0\leq\nu\leq \frac{1}{2}} g_b(\nu)= g_b(\nu_b) 
&=& \frac{2(\sqrt{b}-1)\left(\frac{1}{\log b}-\frac{1}{2(\sqrt{b}-1)}\right)+1}{b^{\frac{1}{\log b}}b^{\frac{-1}{2(\sqrt{b}-1)}}}\\
&=&\frac{\frac{2(\sqrt{b}-1)}{\log b}}{eb^{\frac{-1}{2(\sqrt{b}-1)}}}
=\frac{\left(\sqrt{b}\right)^{\frac{1}{\sqrt{b}-1}}}{e\log{\left(\sqrt{b}\right)}^{\frac{1}{\sqrt{b}-1}}}=S(\sqrt{b}).
\end{eqnarray*}
Thus we have the following inequality.
\begin{equation}
S(\sqrt{b})b^{\nu} \geq \nu b+(1-\nu)-\nu(\sqrt{b}-1)^2.
\end{equation}
In the case of $b=1$, we have the equality in the above inequality, since we have $S(1)=1$.
Replacing $b$ by $\frac{b}{a}$ and then multiplying $a$ to both sides, we have
\begin{equation}
S\left(\sqrt{\frac{a}{b}}\right)a^{1-\nu}b^{\nu} \geq (1-\nu) a+ \nu b-\nu(\sqrt{a}-\sqrt{b})^2,
\end{equation}
since we have $S(x)=S(1/x)$ for $x>0$.

\item[(ii)] For the case of $\nu \geq 1/2$, $r=1-\nu$.
We consider the following function.
$$
h_a(\nu) \equiv \frac{\nu+(1-\nu)a-(1-\nu)(1-\sqrt{a})^2}{a^{1-\nu}},\quad \left(\frac{1}{2}\leq \nu \leq 1 \right).
$$
By the similar way to (i), we have
$$
h_a'(\nu)=0 \Leftrightarrow  \nu =  1-  \left(\frac{1}{\log a}-\frac{1}{2(\sqrt{a}-1)} \right) \equiv \nu_a.
$$
We also find  $\nu_a \in [\frac{1}{2},1]$ and $h'_a(\nu) >0$ for $\nu < \nu_a$ and $h'_a(\nu) <0$ for $\nu > \nu_a$.
Thus the function $h_a(\nu)$ takes a maximum value when  $\nu=\nu_a,\, (a\neq 1)$ and it is calculated by 
$$\max_{\frac{1}{2}\leq \nu \leq 1}h_a(\nu)=h_a(\nu_a)=S(\sqrt{a}).$$ 
Therefore we have the following inequality.
\begin{equation}
S(\sqrt{a})a^{1-\nu} \geq \nu +(1-\nu)a-(1-\nu)(1-\sqrt{a})^2.
\end{equation}
In the case of $a=1$, we have the equality in the above inequality, since we have $S(1)=1$.
Replacing $a$ by $\frac{a}{b}$ and then multiplying $b$ to both sides, we have
\begin{equation}
S\left(\sqrt{\frac{a}{b}}\right)a^{1-\nu}b^{\nu} \geq (1-\nu) a+ \nu b -(1-\nu)(\sqrt{a}-\sqrt{b})^2.
\end{equation}
\end{itemize}
From (i) and (ii), the proof is completed.

\hfill \qed

\begin{Rem}
We easily find that both sides in the inequality (\ref{ineq22}) is less than or equal to thoes  in the inequality (\ref{ineq21}) so that
neither the inequality (\ref{ineq22}) nor the inequality (\ref{ineq21}) is uniformly better than the other.

In addition, our next interest is the ordering between $S\left(\sqrt{\frac{a}{b}}\right) a^{1-\nu}b^{\nu}$ and $(1-\nu) a+ \nu b$.
However we have no ordering between them, because we have the following examples.
For example, let $a=1$ and $b=10$.
If $\nu=0.9$, then $(1-\nu) a+ \nu b - S\left(\sqrt{\frac{a}{b}}\right) a^{1-\nu}b^{\nu}  \simeq -0.246929$.
And if $\nu=0.6$, then $(1- \nu) a+ \nu b - S\left(\sqrt{\frac{a}{b}}\right) a^{1-\nu}b^{\nu}  \simeq  1.71544$.
\end{Rem}

Applying Lemma \ref{lem21}, we have the reverse ratio inequality of the refined Young inequality for positive operators.
\begin{The}
We suppose two invertible positive operators $A$ and $B$ satisfy $0<mI \leq A,B \leq MI$, where $I$ represents an identity operator and $m,M \in \mathbb{R}$. 
For any $\nu\in[0,1]$, we then have
\begin{equation} \label{ineq_reverse_ratio_Young_op}
S(\sqrt{h})A\sharp_{\nu}B \geq (1- \nu) A + \nu B-2r \left(\frac{A+B}{2}-A\sharp_{1/2}B\right),
\end{equation} 
where $h \equiv \frac{M}{m}>1$ and $r \equiv \min\left\{\nu,1-\nu\right\}$.
\end{The}

{\it Proof}:
In Lemma \ref{lem21}, we put $a=1$, then we have for all $b>0$,
$$
S(\sqrt{b})b^{\nu} \geq \nu b +(1-\nu)-r(\sqrt{b}-1)^2
$$
We consider the invertible positive operator $T$ such that $0<mI\leq T\leq MI$.
Then we have the following inequality
\begin{equation}
\max_{m\leq t\leq M} S(\sqrt{t}) T^{\nu} \geq  \nu T + (1-\nu) -r (T-2T^{1/2}+1),
\end{equation}
for any $\nu \in [0,1]$. 
We put $T= A^{-1/2}BA^{-1/2}.$ Since we then have $\frac{1}{h} =\frac{m}{M} \leq A^{-1/2}BA^{-1/2}  \leq \frac{M}{m}=h$,
 we have
\begin{eqnarray*}
&& \max_{\frac{1}{h}\leq t\leq h} S(\sqrt{t}) \left(A^{-1/2}BA^{-1/2} \right)^{\nu} \\
&&\hspace*{-5mm}  \geq  \nu A^{-1/2}BA^{-1/2}   + (1-\nu) 
-r \left\{A^{-1/2}BA^{-1/2}  -2 \left(A^{-1/2}BA^{-1/2}\right)^{1/2}+1\right\}.
\end{eqnarray*}
Note that $h>1$ and $S(x)$ is monotone decreasing for $0<x<1$ and
monotone increasing for $x>1$ \cite{Tom1}.
Thus we have
\begin{eqnarray*}
&& S(\sqrt{h}) \left(A^{-1/2}BA^{-1/2}\right)^{\nu} \\
&& \hspace*{-5mm} \geq \nu A^{-1/2}BA^{-1/2}   + (1-\nu) -r \left\{A^{-1/2}BA^{-1/2}  -2 \left(A^{-1/2}BA^{-1/2}\right)^{1/2}+1\right\}.
\end{eqnarray*}

Multiplying $A^{1/2}$ to the above inequality from both sides, we have the present theorem.

\hfill \qed

\section{A reverse difference inequality for a refined Young inequality}
For the classical Young inequality, the following reverse inequality is known.
For positive real numbers $a,b$ and $\nu\in[0,1]$, M.Tominaga showed the following inequality \cite{Tom1}:
\begin{equation}\label{ineq100} 
L(a,b)\log S\left(\frac{a}{b}\right) \geq (1- \nu) a + \nu b -a^{1-\nu} b^{\nu},
\end{equation}
which is called the converse difference inequality for the Young inequality in \cite{Tom1}

In this section, we show the reverse difference inequality of  the refined Young inequality (\ref{ineq02_rev}).
\begin{Lem}  \label{the100}
For positive real numbers $a, b$ and $\nu \in [0,1]$, we have
\begin{equation}\label{ineq101}
\omega  L(\sqrt{a},\sqrt{b})\log S\left(  \sqrt{\frac{a}{b}}  \right) \geq (1-\nu) a +  \nu b -a^{1-\nu} b^{\nu} -r\left(\sqrt{a}-\sqrt{b}\right)^2,
\end{equation}
where $\omega \equiv \max \left\{\sqrt{a},\sqrt{b}\right\}$.
\end{Lem}
{\it Proof}:

\begin{itemize}
\item[(i)] For the case of $\nu \leq1/2$, $r=\nu$.
We consider the following function.
$$
g_b(\nu) \equiv \nu b +(1-\nu)-b^{\nu}-\nu(\sqrt{b}-1)^2,\quad \left(0\leq \nu \leq \frac{1}{2}\right).
$$
From $g_b'(\nu)=2(\sqrt{b}-1)-b^{\nu}\log b$, we have
$$
g_b'(\nu)=0  \Leftrightarrow \nu = \frac{\log \frac{\sqrt{b}-1}{\log \sqrt{a}}}{\log b} \equiv \nu_b.
$$
We also find that $\nu_b \in [0,\frac{1}{2}]$ by elementaly calculations with the following inequalities:
$$
 1-\frac{1}{\sqrt{b}}\leq \log \sqrt{b} \leq \sqrt{b} -1, \quad (b > 0).
$$
In addition, we have $g_b''(\nu) =-b^{\nu}\left(\log b\right)^2 <0$.
Therefore $g_b$ takes a maximum value when $\nu = \nu_b$, and it is
calculated as $g_b(\nu_b)=L(1,\sqrt{b}) \log S(\sqrt{b})$ by simple but slightly complicated calculations.
Thus we have
$$
L(1,\sqrt{b})\log S\left(\sqrt{b}\right) \geq \nu b+ (1-\nu) -b^{\nu}-\nu(\sqrt{b}-1)^2.
$$
We put $\frac{b}{a}$ instead of $b$ in the above inequality, and then multiplying $a$ to both sides, we have
\begin{equation}  \label{ineq11}
\sqrt{a} L(\sqrt{a},\sqrt{b})\log S\left(\sqrt{\frac{a}{b}}\right) \geq (1-\nu) a + \nu b -a^{1-\nu}b^{\nu}-\nu(\sqrt{a}-\sqrt{b})^2,
\end{equation}
since $L(x,y)=L(y,x)$ and $S(x)=S(1/x)$ for $x>0$.

\item[(ii)] For the case of $\nu \geq1/2$, $r=1-\nu$.
We consider the following function.
$$
h_a(\nu) \equiv \nu  +(1-\nu)a -a^{1-\nu}-(1-\nu)(1-\sqrt{a})^2,\quad \left(\frac{1}{2}\leq \nu \leq 1\right).
$$
By the similar way to (i), we have
$$
h_a'(\nu)=0 \Leftrightarrow  \nu =  1- \frac{\log \frac{\sqrt{a}-1}{\log \sqrt{a}}}{\log a} \equiv \nu_a.
$$
By the similar way to (i), we have $\nu_a \in[\frac{1}{2},1]$ and $h_a''(\nu)=-a^{1-\nu}\left(\log a\right)^2 <0$ so that 
$h_a$ takes a maximum value when $\nu=\nu_a$, and it is calculated as $h_a(\nu_a)  =L(1,\sqrt{a})\log S(\sqrt{a})$.
Thus we have 
$$
L(1,\sqrt{a})\log S(\sqrt{a}) \geq \nu +(1-\nu)a-a^{1-\nu}-(1-\nu)(1-\sqrt{a})^2,
$$
which implies
\begin{equation} \label{ineq12}
\sqrt{b} L(\sqrt{b},\sqrt{a})\log S\left(\sqrt{\frac{a}{b}}\right) \geq (1-\nu) a +  \nu b -a^{1-\nu}b^{\nu}-(1-\nu)(\sqrt{a}-\sqrt{b})^2,
\end{equation}
by replacing $a$ by $\frac{a}{b}$ and then multiplying $b$ to both sides.
\end{itemize}
From the inequalities (\ref{ineq11}) and (\ref{ineq12}), we have the present theorem,
since $L(x,y)=L(y,x)$ and $S(x)=S(1/x)$ for $x>0$.

\hfill \qed

\begin{Rem}   \label{remark_comparison}
We easily find that the right hand side of the inequality (\ref{ineq100}) is greater than that of  the inequality (\ref{ineq101}).
Therefore, if the left hand side of  the inequality (\ref{ineq101}) is  greater than that of  the inequality (\ref{ineq100}),
then Theorem \ref{the100} is trivial one. However, we have not yet found any counter-example such that
\begin{equation} \label{ineq102}
 L(a,b)\log S\left(\frac{a}{b}\right)  \geq \omega L(\sqrt{a},\sqrt{b} )\log S\left( \sqrt{\frac{a}{b}}  \right),
\end{equation}
where $\omega = \max\left\{\sqrt{a},\sqrt{b}\right\}$
for any $a,b >0$ by the computer calculations. Here we give a remark that we have the following inequalities:
\begin{equation}
L(a,b) \leq \omega L(\sqrt{a},\sqrt{b}),\quad and \quad \log S\left(\frac{a}{b}\right) \geq \log S\left(  \sqrt{\frac{a}{b}}  \right)
\end{equation}
for any $a,b >0$. 
At least, we actually have many examples satisfying the inequality (\ref{ineq102}) 
so that we claim that Theorem \ref{the100} is nontrivial as a refinement of the inequality (\ref{ineq100}).

In addition, it is remarkable that we have no ordering between
$$ L(a,b)\log S\left(\frac{a}{b}\right)$$
and
$$\omega L(\sqrt{a},\sqrt{b} )\log S\left(   \sqrt{\frac{a}{b}}  \right)+r\left( \sqrt{a}-\sqrt{b} \right)^2 $$
for any $a,b>0$ and $\nu \in [0,1]$.
Therefore we may claim that Theorem \ref{the100} is also nontrivial from the sense of finding a tighter upper bound of
$(1-\nu) a +  \nu b -a^{1-\nu} b^{\nu} $.

\end{Rem}


Finally we prove the following theorem. It can be proven by the similar method in \cite{Tom1}.

\begin{The}\label{the_reverse_Young_op}
We suppose two invertible positive operators satisfy $0<mI \leq A,B \leq MI$, where $I$ represents an identity operator and $m,M \in \mathbb{R}$. 
For any $\nu\in[0,1]$, we then have
\begin{equation} \label{ineq_reverse_Young_op}
h \sqrt{M} L(\sqrt{M},\sqrt{m})\log S(\sqrt{h}) \geq (1-\nu) A + \nu B-A\sharp_{\nu}B -2r\left(\frac{A+B}{2} -A\sharp_{1/2}B\right),
\end{equation} 
where $h \equiv \frac{M}{m}>1$ and $r \equiv \min\left\{\nu,1-\nu\right\}$.
\end{The}
{\it Proof}:
From the inequality (\ref{ineq101}), we have
\begin{equation}
\omega L(\sqrt{b},1)\log S(\sqrt{b}) \geq \nu b + (1-\nu) -b^{\nu} -r (b-2\sqrt{b}+1),
\end{equation}
for all $\nu\in[0,1]$, putting $b=1$.
We consider the invertible positive operator $T$ such that $0<mI\leq T\leq MI$.
Then we have the following inequality
\begin{equation}
\max_{m\leq t\leq M} \max\left\{\sqrt{t},1\right\} L(\sqrt{t},1)\log S(\sqrt{t}) \geq  \nu T + (1-\nu) -T^{\nu} -r (T-2T^{1/2}+1),
\end{equation}
for any $\nu \in [0,1]$. 
We put $T= A^{-1/2}BA^{-1/2}.$ Since we then have $\frac{1}{h} =\frac{m}{M} \leq A^{-1/2}BA^{-1/2}  \leq \frac{M}{m}=h$,
 we have
\begin{eqnarray*}
&& \max_{\frac{1}{h}\leq t\leq h} \max\left\{\sqrt{t},1\right\} L(\sqrt{t},1)\log S(\sqrt{t}) \\
&&\hspace*{-5mm}  \geq  \nu A^{-1/2}BA^{-1/2}   + (1-\nu) 
-\left(A^{-1/2}BA^{-1/2}\right)^{\nu} -r \left\{A^{-1/2}BA^{-1/2}  -2 \left(A^{-1/2}BA^{-1/2}\right)^{1/2}+1\right\}.
\end{eqnarray*}
Note that $h>1$ and $L(u,1)$ is monotone increasing function for $u>0$. In addition, we note that $S(x)$ is monotone decreasing for $0<x<1$ and
monotone increasing for $x>1$ \cite{Tom1}.
Thus we have
\begin{eqnarray*}
&& \sqrt{h}L(\sqrt{h},1)\log S(\sqrt{h}) \\
&& \hspace*{-5mm} \geq \nu A^{-1/2}BA^{-1/2}   + (1-\nu) -\left(A^{-1/2}BA^{-1/2}\right)^{\nu} -r \left\{A^{-1/2}BA^{-1/2}  -2 \left(A^{-1/2}BA^{-1/2}\right)^{1/2}+1\right\}.
\end{eqnarray*}

Multiplying $A^{1/2}$ to the above inequality from both sides, we have
$$
\sqrt{h} L(\sqrt{h},1)\log S(\sqrt{h})A \geq (1-\nu) A+ \nu B- A\sharp_{\nu}B -2r\left(\frac{A+B}{2}-A\sharp_{1/2}B\right).
$$
Since the left hand side in the above inequality is less than 
$$\sqrt{h} L(\sqrt{h},1)\log S(\sqrt{h})M= h \sqrt{M}L(\sqrt{M},\sqrt{m})\log S(\sqrt{h}),$$
the proof is completed.

\hfill \qed

\begin{Rem}

As mentioned in Remark \ref{remark_comparison}, we have not yet found the ordering between
the right hand side of the inequality (\ref{ineq_reverse_Young_op}) and that of the inequality (\ref{prop_rev_diff_ineq}).
Therefore Theorem \ref{the_reverse_Young_op} is not a trivial result.
\end{Rem}


\section{Concluding remarks}
As we have seen, we gave refined Young inequalities for two positive operators.
In addition, we gave reverse ratio type inequalities and reverse difference type inequalities for the refined Young inequality for positive operators.
Closing this paper, we shall give a refinement of the weighted arithmetic-geometric mean inequality for $n$ real numbers by a simple proof.
\begin{Prop}  \label{the_gen01}
Let  $a_1,\cdots,a_n\geq 0$ and $p_1,\cdots,p_n > 0$ with $\sum_{j=1}^n p_j=1$ and  $\lambda  \equiv \min \left\{ p_1,\cdots ,p_n \right\}$.
If we assume that the multiplicity attaining $\lambda$ is $1$,
then we have
\begin{equation} \label{gen_ineq01}
\sum_{i=1}^n p_i a_i  -   \prod_{i=1}^n a_i^{p_i}  \geq n \lambda  \left(\frac{1}{n}\sum_{i=1}^n a_i-  \prod_{i=1}^n a_i^{1/n}     \right),
\end{equation}
with equality if and only if $a_1=\cdots =a_n$.
\end{Prop}

{\it Proof}:
We suppose $\lambda = p_j$.  For any $j=1,\cdots,n$, 
we then have
\begin{eqnarray*}
 \sum_{i=1}^n p_i a_i  - p_j \left(\sum_{i=1}^n a_i - n \prod_{i=1}^n a_i^{1/n}  \right) 
&=& np_j \left( \prod_{i=1}^n a_i^{1/n}      \right)  + \sum_{i=1,i\neq j}^n ( p_i-p_j) a_i \\
&\geq&  \prod_{i=1,i\neq j}^n \left(   a_1^{1/n} \cdots a_n^{1/n}  \right)^{np_j}  a_i^{p_i-p_j} \\
&=& a_1^{p_1}  \cdots a_n^{p_n}.
\end{eqnarray*}
In the above process, the classical weighted arithmetic-geometric mean inequality \cite{HLP,Bu} for
$a_1,\cdots,a_n\geq 0$ and  $p_1,\cdots,p_n > 0$ with  $\sum_{j=1}^n p_j=1$;
\begin{equation} \label{classical_w_ineq}
\sum_{j=1}^n p_j a_j \geq \prod_{j=1}^n a_j^{p_j},
\end{equation}
with equality if and only if $a_1=\cdots =a_n$, was used.
We note that $p_i-p_j>0$ from the assumption of the proposition. 
The equality in the inequality (\ref{gen_ineq01}) holds if and only if 
$$(a_1a_2\cdots a_n)^{\frac{1}{n}}=a_1=a_2=\cdots=a_{j-1}=a_{j+1}=\cdots=a_n$$ 
by the equality condition of the classical weighted arithmetic-geometric mean inequality  (\ref{classical_w_ineq}). 
Therefore $a_1=a_2=\cdots=a_{j-1}=a_{j+1}=\cdots=a_n\equiv a$, then we have $a_j^{\frac{1}{n}}a^{\frac{n-1}{n}}=a$ from the first equality.
Thus we have $a_j=a$, which completes the proof. 
\hfill \qed

The inequality (\ref{gen_ineq01}) gives a refinement of the classical weighted arithmetic-geometric mean inequality (\ref{classical_w_ineq}). 
At the same time, it gives a natural generalization of the inequality (\ref{ineq02_rev}) proved in \cite{KM}.
It is also notable that Proposition \ref{the_gen01} can be proven by using the bounds for
the normalized Jensen functional, which were obtained by S.S.Dragomir in \cite{Dra}, except for the equality condtions.
Note that the inequality (\ref{gen_ineq01}) itself holds without the assumption that the multiplicity attaining $\lambda$ is $1$.
In addition, when we do not impose on this assumption, the equality in the inequality (\ref{gen_ineq01}) holds if $p_i=\frac{1}{n}$ for all $i=1,\cdots,n$.

\section*{Ackowledgement}
The author thanks Professor S. M. Sitnik for letting me know the valuable infromation on the reference \cite{BK} and 
Dr.F.C.Mitroi for giving me the valuable comment on the equality condition of Propisition \ref{the_gen01}.
The author also thanks the anonymous referee for valuable comments to improve the manuscript.
The author was supported in part by the Japanese Ministry of Education, Science, Sports and Culture, Grant-in-Aid for 
Encouragement of Young Scientists (B), 20740067.


\begin{thebibliography}{99}
\bibitem{BK} N.A.Bobylev and M. A. Krasnoselsky, Extremum Analysis (degenerate cases), Moscow, preprint, 1981, 52 pages, (in Russian).
\bibitem{KM} F.Kittaneh and Y.Manasrah, Improved Young and Heinz inequalities for matrices, J.Math.Anal.Appl.,Vol.36(2010), pp.262-269.
\bibitem{IN} S.Izumino and N.Nakamura, Weighted Geometric Means of Positive Operators, Kyungpook Math. J. 50(2010), pp.213-228.
\bibitem{Tom1} M.Tominaga, Specht's ratio in the Young inequality, Sci.Math.Japon.,Vol.55(2002),pp.583-588.
\bibitem{Specht} W.Specht, Zer Theorie der elementaren Mittel, Math.Z.,Vol.74(1960),pp.91-98.
\bibitem{JIFUJII} J.I.Fujii, S.Izumino and Y.Seo, Determinant for positive operators and Specht's theorem, Sci.Math.Japon.,Vol.1(1998),pp.307-310.
\bibitem{KA} F.Kubo and T.Ando, Means of positive operators, Math. Ann.,Vol.264(1980),pp.205-224.
\bibitem{FY} T.Furuta and M.Yanagida, Generalized means and convexity of inversion for positive operators, Amer.Math.Monthly, Vol.105 (1998),pp.258-259.
\bibitem{Furuta} T.Furuta, Invitation to linear operators: From matrix to bounded linear operators on a Hilbert space, Taylor and  Francis, 2002.
\bibitem{Furuta1} T.Furuta, The H\"older-McCarthy and the Young inequalities are equivalent for Hilbert space operators, Amer.Math.Monthly, Vol.108(2001),pp.68-69.
\bibitem{Mc} C.A.McCarthy, $C_p$, Israel J. Math.,Vol.5(1967),pp.249-271.
\bibitem{HLP} G.H.Hardy, J.E.Littlewood and G.P\'olya, Inequalities, Cambridge University Press, 1952.
\bibitem{Bu} P.S.Bullen, Handbook of Means and Their Inequalities, Kluwer Acad. Publ., Dordrecht, 2003. 
\bibitem{Dra} S.S.Dragomir, Bounds for the normalised Jensen functional, Bull. Australian Math. Soc., Vol.74(2006), pp.471-478.    

\end{thebibliography}
\end{document}